\documentclass[12pt]{article}

\usepackage{a4,graphicx,amsmath,amsfonts,amssymb,mathrsfs,bbm}

\newcommand{\real}{\mathbbm{R}}
\newcommand{\complex}{\mathbbm{C}}
\newcommand{\nat}{\mathbbm{N}}
\newcommand{\ltwo}{\mathscr{L}^2(\Pi,\rho)}
\newcommand{\linf}{\mathscr{L}^{\infty}(\Pi,\rho)}
\newcommand{\cl}{\mathscr{C}^{\ell}(\Pi)}
\newcommand{\czero}{\mathscr{C}^{0}(\Pi)}
\newcommand{\innerprod}[2]{\left\langle #1, #2 \right\rangle}

\newtheorem{definition}{Definition}
\newtheorem{theorem}{Theorem}

\newtheorem{lemma}{Lemma}

\begin{document}
\begin{center}
  {\bf \Large Stability preservation in \\[0.5ex]
    stochastic Galerkin projections \\[1.5ex]
    of dynamical systems}
\vspace{10mm}

{\large Roland~Pulch} \\[1ex]
{\small Institute of Mathematics and Computer Science,
  University of Greifswald, \\
Walther-Rathenau-Str.~47, 17489 Greifswald, Germany. \\
Email: {\tt pulchr@uni-greifswald.de}}

\vspace{5mm}

{\large Florian Augustin} \\[1ex]
{\small Massachusetts Institute of Technology \\
Cambridge, MA 02139, United States. \\
Email: {\tt fmaugust@mit.edu }}

\end{center}

\bigskip\bigskip


\begin{center}
{Abstract}

\begin{tabular}{p{13cm}}
  In uncertainty quantification, critical parameters of mathematical models
  are substituted by random variables.
  We consider dynamical systems composed of ordinary differential equations.
  The unknown solution is expanded into an orthogonal basis of the
  random space, e.g., the polynomial chaos expansions.
  A Galerkin method yields a numerical solution of the stochastic model.
  In the linear case, the Galerkin-projected system may be unstable,
  even though all realizations of the original system are asymptotically
  stable.
  We derive a basis transformation for the state variables in the
  original system,
  which guarantees a stable Galerkin-projected system.
  The transformation matrix is obtained from a symmetric decomposition
  of a solution of a Lyapunov equation.
  In the nonlinear case, we examine stationary solutions
  of the original system.
  Again the basis transformation preserves the asymptotic stability of the
  stationary solutions in the stochastic Galerkin projection.
  We present results of numerical computations for both
  a linear and a nonlinear test example.
  \bigskip

  \bigskip

\end{tabular}
\end{center}

\clearpage


\section{Introduction}
Uncertainty quantification (UQ) examines the dependence of outputs
on vague input parameters in mathematical models, see~\cite{sullivan}.
Often the uncertain parameters are replaced by random variables or
random processes, resulting in a stochastic problem.
We consider dynamical systems consisting of ordinary differential equations
(ODEs) with random parameters.
The state variables can be expanded into a series of orthogonal
basis functions, where often polynomials are applied (polynomial chaos),
see~\cite{augustinetal,ernst,xiu-book}.
Stochastic Galerkin methods or stochastic collocation techniques
yield numerical solutions of unknown coefficient functions.
We focus on the stochastic Galerkin approach,
see~\cite{augustin-rentrop,pulch11,pulch14},
in this paper.

Sonday et al.~\cite{sonday} analyzed the spectrum of a Jacobian matrix
of a Galerkin-projected (nonlinear) system of ODEs.
The results indicate that a preservation of stability is not
guaranteed in the Galerkin projection even if the original system
is asymptotically stable. 
Even though a loss of stability happens rather seldom, this change
in stability leads to unexpected and erroneous results.
Therefore, we derive a technique, which guarantees the asymptotic stability
in the Galerkin-projected system provided that the original system is
asymptotically stable.

Prajna~\cite{prajna} designed an approach to preserve stability in a
pro\-jec\-tion-based model order reduction of a (nonlinear) system of ODEs.
Therein, a dynamical system is reduced to a smaller dynamical system.
We apply a similar strategy in the stochastic Galerkin method,
where a random dynamical system is projected to a larger deterministic
dynamical system.
The stability-preserving technique employs a basis transformation of
the original parameter-dependent system,
where the transformation matrix is derived from the solution of a
Lyapunov equation.

We construct and investigate this transformation for
linear dynamical systems in detail.
A proof of the stability preservation is given for the stabilized stochastic
Galerkin method.
Furthermore, we consider the asymptotic stability of stationary solutions
(equilibria) for autonomous nonlinear dynamical systems.
In~\cite{pulch13}, existence and convergence of stationary solutions
was analyzed in the stochastic Galerkin-projected systems.
The Galerkin system exhibits equilibria, which yield approximations
to the random-dependent equilibria of the original system.
The approximations converge in mean square to the exact equilibria. 
Now we apply the basis transformation to guarantee the stability of
the stationary solutions in the stochastic Galerkin-projected system.

The paper is organized as follows.
The stochastic Galerkin approach is described for the linear case
in Section~\ref{sec:problem-def}.
The stability-preserving projection is derived and analyzed in
Section~\ref{sec:stability}.
An analogous stabilization is specified for the nonlinear case
in Section~\ref{sec:nonlinear}.
Finally, Section~\ref{sec:example} includes numerical results for
both a linear and a nonlinear test example.


\section{Problem definition}
\label{sec:problem-def}
The class of linear problems under investigation is described
in this section.

\subsection{Linear dynamical systems and stability}
We consider a linear dynamical system of the form
\begin{equation} \label{linear-ode}
  \dot{x}(t,p) = A(p) x(t,p) + s(t,p) ,
\end{equation}
where the matrix $A : \Pi \rightarrow \real^{n \times n}$ and the
vector $s : [0,\infty) \times \Pi \rightarrow \real^n$ depend
on parameters $p \in \Pi$ for some subset $\Pi \subseteq \real^q$.
Consequently, the state variables
$x : [0,\infty) \times \Pi \rightarrow \real^n$
are also parameter-dependent.
Initial value problems are determined by
$$ x(0,p) = x_0(p) $$
with a given function $x_0 : \Pi \rightarrow \real^n$.
Since we are investigating stability properties, let, without loss of
generality, $s \equiv 0$ in the system~(\ref{linear-ode}).

To analyze the stability, we recall some general properties of matrices.
\begin{definition} \label{def:stable-matrix}
  Let $A \in \real^{n \times n}$ and $\lambda_1,\ldots,\lambda_n \in \complex$
  be its eigenvalues.
  The {\em spectral abscissa} of the matrix~$A$ reads as
  $$ \alpha (A) := \max
  \left\{ {\rm Re}(\lambda_1) , \ldots , {\rm Re} (\lambda_n) \right\} . $$
  $A$ is called a {\em stable matrix}, if it holds that $\alpha (A) < 0$. 
\end{definition}
A linear dynamical system $\dot{x} = Ax$ is asymptotically stable
if and only if the included matrix~$A$ is stable.
We assume that the matrices~$A(p)$ in the system~(\ref{linear-ode})
are stable for all~$p \in \Pi$ in the following.

\subsection{Stochastic modeling and orthogonal expansions}
Now we assume that the parameters in equation (\ref{linear-ode}) are affected by uncertainties.
In uncertainty quantification, 
the parameters are replaced by independent random variables
$p : \Omega \rightarrow \Pi$ on some
probability space $(\Omega,\mathscr{A},\mu)$.
Let a joint probability density function $\rho : \Pi \rightarrow \real$
be given.
Without loss of generality, we assume $\Pi = {\rm supp}(\rho)$,
because the parameter space~$\Pi$ can be restricted to the support
of~$\rho$ otherwise.
For a measurable function $f: \Pi \rightarrow \real$,
the expected value reads as
\begin{equation} \label{expectedvalue}
  \mathbb{E} \left[ f \right] :=
  \int_{\Omega} f(p(\omega)) \; \mbox{d}\mu(\omega) =
  \int_{\Pi} f(p) \rho(p) \; \mbox{d}p
\end{equation}
provided that the integral exists.
The Hilbert space
$$ \ltwo :=
\left\{ f : \Pi \rightarrow \real \; : \;
f \; \mbox{measurable and} \; \mathbb{E} \left[ f^2 \right] < \infty
\right\} $$
is equipped with the inner product
$$ \innerprod{f}{g} \; = \int_{\Pi} f(p) g(p) \rho(p) \; \mbox{d}p
\qquad \mbox{for} \;\; f,g \in \ltwo . $$
Let a complete orthonormal system $(\Phi_i)_{i \in \nat}$ be given.
Thus the basis functions $\Phi_i : \Pi \rightarrow \real$ satisfy
\begin{equation} \label{orthogonal}
  \innerprod{\Phi_i}{\Phi_j} \; = \left\{
  \begin{array}{ll}
    0 & \mbox{for} \;\; i \neq j , \\
    1 & \mbox{for} \;\; i = j . \\
  \end{array}
  \right.
\end{equation}
Assuming $x_k(t,\cdot) \in \ltwo$ for each component
$k = 1,\ldots,n$ and each time point~$t$,
the state variables of the system~(\ref{linear-ode})
can be expanded into a series
\begin{equation} \label{expansion}
  x(t,p) = \sum_{i=1}^{\infty} v_i(t) \Phi_i(p) .
\end{equation}
The coefficient functions $v_i : [0,\infty) \rightarrow \real^n$
are defined by
\begin{equation} \label{coefficients}
  v_{i,k} = \; \innerprod{x_k(t,\cdot)}{\Phi_i(\cdot)}
  \qquad \mbox{for} \;\; k=1,\ldots,n .
\end{equation}
The series~(\ref{expansion}) converges in the norm of
$\ltwo$ point-wise for each~$t$.
Often polynomials are used as basis functions following the
concepts of the generalized polynomial chaos (gPC).
More details can be found in~\cite{xiu-book}.

\subsection{Galerkin projection of linear dynamical systems}
Stochastic Galerkin methods and stochastic collocation techniques yield
approximations of the coefficient functions~(\ref{coefficients})
in the expansion~(\ref{expansion}),
see~\cite{augustin-rentrop,pulch11,pulch14,pulch-maten-augustin}.
We apply the stochastic Galerkin approach,
where Equation~(\ref{linear-ode}) is projected onto a finite subset $\{ \Phi_1 , \ldots , \Phi_m \}$ of basis functions. The stochastic process~(\ref{expansion}) is approximated by a
truncated expansion
$$ \hat{x}^{(m)} (t,p) =
\sum_{i=1}^{m} \hat{v}_i(t) \Phi_i(p) . $$
The Galerkin projection of the dynamical system~(\ref{linear-ode}), neglecting the term~$s$, results in the larger linear dynamical system
\begin{equation} \label{system-galerkin}
  \dot{\hat{v}}(t) = \hat{A} \hat{v}(t) ,
\end{equation}
whose solution $\hat{v} = ( \hat{v}_1^\top , \ldots , \hat{v}_m^\top )^\top$
represents an approximation of the exact coefficient
functions~(\ref{coefficients}).
The matrix $\hat{A} \in \real^{mn \times mn}$ is defined by its
minors $\hat{A}_{ij} \in \real^{n \times n}$ with
$$ \hat{A}_{ij} = \mathbb{E} \left[ A \Phi_i \Phi_j \right]
\qquad \mbox{for} \;\; i,j = 1,\ldots,m $$
using the matrix~$A$ from~(\ref{linear-ode}).
Therein, the expected value, see~(\ref{expectedvalue}), is applied
com\-ponent\-wise.
If the matrix $A(p)$ is symmetric for almost all~$p$,
then the matrix $\hat{A}$ is also symmetric.
Otherwise, the matrix $\hat{A}$ is unsymmetric, which is the case in many situations.

The convergence properties of the stochastic Galerkin approach
are not investigated in this paper.
Alternatively, we examine the stability properties.
The analysis in~\cite{sonday} shows that the matrix~$\hat{A}$ may be
unstable even though $A(p)$ is stable for strictly all~$p$
with respect to Definition~\ref{def:stable-matrix}.
Yet the stability is guaranteed in the case of normal matrices $A(p)$
for almost all~$p$.
Even though stability can be lost for non-normal matrices,
the stochastic Galerkin method is still convergent on compact time intervals
under usual assumptions.

\subsection{Basis transformations}
We consider a transformation of the linear dynamical
system~(\ref{linear-ode}), with $s \equiv 0$, to an equivalent system
\begin{equation} \label{system-transformed}
  \dot{y}(t,p) = B(p) y(t,p)
\end{equation}
with $y(t,p) := T(p) x(t,p)$ and transformation matrices
$T : \Pi \rightarrow \real^{n \times n}$ being point-wise non-singular.
It holds that
\begin{equation} \label{similarity-trafo}
  B(p) = T(p) A(p) T(p)^{-1}
\end{equation}
for each $p \in \Pi$.
The operation~(\ref{similarity-trafo}) represents a similarity
transformation, i.e., the spectra of the matrices $A(p)$ and $B(p)$ 
coincide.

If the stochastic Galerkin system~(\ref{system-galerkin}) is unstable,
then our aim is to identify a basis transformation given by a matrix
$T : \Pi \rightarrow \real^{n \times n}$
such that the Galerkin projection of the
dynamical system~(\ref{system-transformed}) yields a stable system.
The following properties of the basis transformation are required:
\begin{enumerate}
\item 
  $T$ has to be non-constant in the variable~$p$.
  The Galerkin approach is invariant with respect to constant
  basis transformations and thus
  the stability properties cannot be changed.
\item
  If $A \in \mathscr{C}^{\ell}(\Pi)^{n \times n}$,
  then $T \in \mathscr{C}^{\ell}(\Pi)^{n \times n}$
  is required to guarantee $B \in \mathscr{C}^{\ell}(\Pi)^{n \times n}$
  and thus $y(t,\cdot) \in \mathscr{C}^{\ell}(\Pi)^n$ for each~$t$. 
  The convergence rate of orthogonal (gPC) expansions depends on
  the order of differentiability in the random-dependent functions, see
  \cite[p.~154]{augustinetal} and \cite[p.~33]{xiu-book}.
\end{enumerate}


\section{Stability preservation}
\label{sec:stability}
We derive a concept to guarantee the stability of the dynamical system
obtained by the stochastic Galerkin projection.

\subsection{General results}
In this subsection, a constant matrix~$A \in \real^{n \times n}$ is considered.
If the matrix~$A$ is unsymmetric, then the following definition
allows for further investigations.
\begin{definition} \label{def:symm-part}
  The {\em symmetric part} of a matrix $A \in \real^{n \times n}$ reads as
  $$ A_{\rm sym} := \textstyle \frac{1}{2} (A + A^\top) . $$ 
\end{definition}
The symmetric part of~$A$ is negative definite if and only if $A + A^\top$
is negative definite.
We will apply the following well-known property later.
\begin{lemma} \label{lemma:stable-matrix}
  If the symmetric part of $A \in \real^{n \times n}$ is negative
  definite, then $A$ is a stable matrix.
\end{lemma}

\underline{Proof:} \nopagebreak

The spectral abscissa is bounded by $\alpha(A) \le \mu (A)$ for an
arbitrary logarithmic norm~$\mu$.
The logarithmic norm associated with the Euclidean vector norm reads as
$\mu(A) = \alpha (A_{\rm sym})$, 
see~\cite[p.~61]{hairer1}. 
The negative definiteness of the symmetric part implies
$\alpha (A_{\rm sym})<0$.
It follows that $\alpha(A) < 0$.
\hfill $\Box$

\medskip

Assuming that the symmetric part of a given matrix~$A$ is not
negative definite, we construct a transformation matrix $T$ such that
the symmetric part of the similarity-transformed matrix $B = TAT^{-1}$
becomes negative definite.

\clearpage

\begin{theorem} \label{thm:symm-part}
  Let $A \in \real^{n \times n}$ be a stable matrix and $Q \in \real^{n \times n}$
  be a symmetric positive definite matrix.
  The Lyapunov equation
  \begin{equation} \label{lyapunov}
    A^\top M + M A + Q = 0
  \end{equation}
  has a unique symmetric positive definite solution $M \in \real^{n \times n}$.
  For a symmetric decomposition $M = L L^\top$ with $L \in \real^{n \times n}$,
  a similarity transformation yields the matrix
  \begin{equation} \label{B-trafo}
    B := L^\top A L^{-\top} ,
  \end{equation}
  which features a negative definite symmetric part.
\end{theorem}

\underline{Proof:} \nopagebreak

The stability of the matrix~$A$ guarantees existence and uniqueness
of a solution~$M$ for the Lyapunov equation~(\ref{lyapunov}),
see~\cite[p.~303]{hammarling}. 
The symmetric matrix~$M$ is positive definite,
because $Q$ is positive definite.
The symmetric part of the matrix~(\ref{B-trafo}) becomes
(neglecting the factor~$\frac{1}{2}$)
$$ \begin{array}{rcl}
  B + B^\top & = & L^\top A L^{-\top} + L^{-1} A^\top L \\[1ex]
  & = & L^{-1} ( M A + A^\top M ) L^{-\top} \\[1ex]
  & = & - L^{-1} Q L^{-\top} . \\
\end{array} $$
The matrix $- L^{-1} Q L^{-\top}$ is negative definite due to the
positive-definiteness of $Q$, since
$$ z^\top (- L^{-1} Q L^{-\top}) z = - (L^{-\top} z)^\top Q (L^{-\top} z) < 0 $$
for all $z \in \real^n\backslash\{0\}$.
\hfill $\Box$

\medskip

The proof of Theorem~\ref{thm:symm-part} follows mainly the steps
in~\cite[Thm.~5]{prajna}.
However, a symmetric decomposition $M = M^{\frac{1}{2}} M^{\frac{1}{2}}$
is assumed in~\cite{prajna}, which requires the computation of
all eigenvalues and eigenvectors.
Furthermore, this decomposition is not unique in the case of
multiple eigenvalues.
Theorem~\ref{thm:symm-part} holds true for arbitrary
symmetric decompositions of~$M$. 
In particular, we may apply the Cholesky factorization
$M = L L^\top$, where $L$ becomes a unique lower triangular matrix with
strictly positive diagonal elements,
see~\cite[p.~204]{stoerbulirsch}. 
Efficient algorithms are available to compute the Cholesky factor
without first finding~$M$, see~\cite{hammarling}.
More details on Lyapunov equations can be found in~\cite{gajic},
for example.

\subsection{Stability-preserving transformation}
\label{sec:galerkin-trafo}
The linear dynamical system~(\ref{linear-ode})
is assumed to involve stable matrices $A(p)$ for all~$p \in \Pi$.
In view of~(\ref{lyapunov}), we use the parameter-dependent
Lyapunov equations
\begin{equation} \label{lyapunov-parameter}
  A(p)^\top M(p) + M(p) A(p) + Q(p) = 0
  \qquad \mbox{for each} \;\; p \in \Pi .
\end{equation}
Let the matrices $Q(p)$ be symmetric positive definite for all~$p$.
Constant choices $Q(p) \equiv Q_0$ are admissible.
Consequently, the system (\ref{lyapunov-parameter}) has a unique symmetric positive definite solution~$M(p)$ for each~$p$.
We require a symmetric decomposition of~$M(p)$ for each~$p$.
Concerning the smoothness, we demonstrate a property of the
Cholesky factorization.
The proof follows the steps in~\cite[p.~295]{schatzman},
where the continuity of this decomposition is shown.

\begin{lemma} \label{lemma:cholesky}
If $A \in \cl^{n \times n}$ and $A(p)$ is symmetric as well as
positive definite for all $p \in \Pi$, then the
Cholesky decomposition $A = L L^\top$ satisfies
$L \in \cl^{n \times n}$.
\end{lemma}

\underline{Proof:} \nopagebreak

We use induction with respect to $n$.
For $n=1$, we obtain $A(p) = (\alpha(p))$ with $\alpha(p) > 0$
for all~$p$.
It follows that $L(p) = (\sqrt{\alpha(p)})$.
Thus $L \in \cl^{1 \times 1}$ is satisfied,
because the square root is differentiable to arbitrary order
for positive real numbers.
Now let the assumption be valid for $n-1$.
We partition a matrix $A(p) \in \mathbbm{R}^{n \times n}$ and its
Cholesky decomposition into
$$ A(p) = \begin{pmatrix}
  \alpha(p) & r(p) \\
  r(p)^\top & \bar{A}(p) \\
\end{pmatrix}
\quad \mbox{and} \quad
 L(p) = \begin{pmatrix}
  \beta(p) & 0 \\
  s(p)^\top & B(p) \\
 \end{pmatrix} . $$
Since $A(p)$ is positive definite, it holds that $\alpha(p) > 0$ for all~$p$.
We obtain
$\beta(p) = \sqrt{\alpha(p)}$,
$s(p) = \frac{r(p)}{\sqrt{\alpha(p)}}$ 
and
$B(p) B(p)^\top = \bar{A}(p) - \frac{r(p)^\top r(p)}{\alpha(p)} =: F(p)$.
The mapping $p \mapsto F(p)$ is in $\cl^{n \times n}$ due to
$A \in \cl^{n \times n}$.
Hence the mapping $p \mapsto B(p)$ is in $\cl^{n \times n}$ by the
assumption in the induction.
Note that the operations used to compute $\beta$ and $s$
are differentiable to arbitrary orders.
\hfill $\Box$

\medskip

The following result guarantees the preservation of smoothness.

\begin{lemma} \label{lemma:smoothness}
  If $A,Q \in \cl^{n \times n}$, 
  then it follows that $B \in \cl^{n \times n}$ for the
  transformed matrix~(\ref{similarity-trafo}) using the
  Cholesky factorization of the solutions from the
  Lyapunov equations~(\ref{lyapunov-parameter}).
\end{lemma}

\underline{Proof:} \nopagebreak

Each Lyapunov equation~(\ref{lyapunov-parameter}) represents a
larger linear system of algebraic equations.
Assuming $A,Q \in \cl^{n \times n}$, Cramer's rule implies that the
entries in the solution~$M$ are also in $\cl$.
Lemma~\ref{lemma:cholesky} yields the differentiability~$L \in \cl^{n \times n}$.
The inverse matrix inherits the smoothness $L^{-1} \in \cl^{n \times n}$,
because it holds that
$L^{-1}(p) = \frac{{\rm adj}(L(p))}{\det(L(p))}$ for each~$p$
with the adjoint matrix.
Now formula~(\ref{similarity-trafo}) demonstrates that
$B \in \cl^{n \times n}$.
\hfill $\Box$

\medskip

In particular, Lemma~\ref{lemma:smoothness} is valid in
the case of continuous functions ($\ell = 0$).
Furthermore, just measurable matrix-valued functions $A,Q$ imply
measurable matrix-valued functions $M,L$.

If the entries of $A$ and $Q$ are polynomials in the variable~$p$,
then the entries of $M$ become rational functions.
However, the entries of the factor~$L$ in a symmetric decomposition
are not rational functions in general,
because a square root is typically applied somewhere.
Consequently, the transformed matrix~(\ref{B-trafo}) does not represent
a rational function in general.

We transform the original dynamical system~(\ref{linear-ode})
into~(\ref{system-transformed}) using the transformation matrix
$T(p) := L(p)^\top$
for an arbitrary symmetric decomposition $M(p) = L(p) L(p)^\top$.
The main result is formulated now.
\begin{theorem} \label{thm:main}
  Let $A,Q : \Pi \rightarrow \real^{n \times n}$ be measurable functions
  with $A(p)$ stable
  and $Q(p)$ symmetric as well as positive definite for almost all
  $p \in \Pi$.
  The Lyapunov equation~(\ref{lyapunov-parameter}) yields a unique
  solution $M(p)$ for almost all~$p$. 
  A symmetric decomposition $M(p) = L(p) L(p)^\top$ is considered
  for almost all~$p$.
  Furthermore, let $A \in \mathscr{L}^{q_1}(\Pi,\rho)^{n \times n}$,
  $L \in \mathscr{L}^{q_2}(\Pi,\rho)^{n \times n}$,
  $L^{-1} \in \mathscr{L}^{q_3}(\Pi,\rho)^{n \times n}$
  with $q_k \in [1,\infty]$ for $k=1,2,3$ and
  $\frac{1}{2} = \frac{1}{q_1} + \frac{1}{q_2} + \frac{1}{q_3}$.
  Using $T(p) := L(p)^\top$, the Galerkin pro\-jec\-tion of the
  transformed system~(\ref{system-transformed})
  with the matrix~(\ref{similarity-trafo})
  produces an asymptotically stable linear dynamical system.
\end{theorem}

\underline{Proof:} \nopagebreak

The generalized H\"older inequality guarantees
$B = TAT^{-1} \in \ltwo^{n \times n}$ due to the
regularity assumptions for $A,L,L^{-1}$.
The stochastic Galerkin approach applied to the
transformed system~(\ref{system-transformed})
including the matrix~(\ref{similarity-trafo}) yields
a dynamical system $\dot{\hat{v}} = \hat{B} \hat{v}$
with the matrix $\hat{B} \in \real^{mn \times mn}$.
The minors $\hat{B}_{ij} \in \real^{n \times n}$ read as
$\hat{B}_{ij} = \mathbb{E} \left[ B \Phi_i \Phi_j \right]$
for $i,j = 1,\ldots,m$.
We investigate the symmetric part of the matrix~$\hat{B}$.
The minors of the symmetric part $\hat{B} + \hat{B}^\top$ become
$$ \hat{B}_{ij} + (\hat{B}_{ji})^\top =
  \mathbb{E} \left[ B \Phi_i \Phi_j \right] +
  \mathbb{E} \left[ B^\top \Phi_j \Phi_i \right]
  = \mathbb{E} \left[ (B + B^\top) \Phi_i \Phi_j \right] . $$
Given the transformation~(\ref{similarity-trafo}) with
$T(p) = L(p)^\top$, Theorem~\ref{thm:symm-part} shows that the
symmetric part, $B(p) + B(p)^\top$, of the matrix~(\ref{B-trafo}) is negative definite for almost all~$p$.
Let $z = (z_1^\top,\ldots,z_m^\top)^\top \in \real^{mn}$
with $z_1,\ldots,z_m \in \real^n$.
It follows that
$$ \begin{array}{rcl}
  z^\top ( \hat{B} + \hat{B}^\top ) z & = & \displaystyle
  \sum_{i,j=1}^m z_i^\top \left( \hat{B}_{ij} + \hat{B}_{ji}^\top \right) z_j
  \;\; = \;\;
  \mathbb{E} \left[ \sum_{i,j=1}^m z_i^\top (B + B^\top) z_j \Phi_i \Phi_j\right] \\
  & = & \displaystyle
  \mathbb{E} \left[ \left( \sum_{i=1}^m z_i \Phi_i \right)^\top
    (B + B^\top) \left( \sum_{j=1}^m z_j \Phi_j \right) \right]
  \;\; \le \;\; 0 . \\
  \end{array} $$ 
Since the basis functions are linearly independent, it holds that
$$ \tilde{z} := \sum_{i=1}^m z_i \Phi_i
\in \ltwo \backslash \{ 0 \}
\qquad \mbox{for} \;\; z \neq 0 . $$
Assuming $z \neq 0$,
the function~$\tilde{z}$ is non-zero on a subset $U \subset \Pi$
with $\mu(U) > 0$ for the probability measure~$\mu$.
It follows that the above expected value becomes strictly negative.
Hence the symmetric part $\hat{B} + \hat{B}^\top$ is negative definite.
Lemma~\ref{lemma:stable-matrix} shows that $\hat{B}$ is a stable matrix.
Consequently, the dynamical system $\dot{\hat{v}} = \hat{B} \hat{v}$
is asymptotically stable.
\hfill $\Box$

\medskip

The above result is independent of the choice of orthogonal basis functions.
Moreover, the system of basis functions is not required to be complete.
Note that any symmetric decomposition can be used satisfying
the suppositions in Theorem~\ref{thm:main}. 

Concerning the regularity assumptions,
the choice $q_1=q_2=q_3=6$ is admissible, for example.
Furthermore, $q_i = \infty$ can be chosen for one particular~$i$,
with $\frac{1}{q_i} = 0$.
It holds that
$\linf \subset \mathscr{L}^q(\Pi,\rho)$
for any $q \in [1,\infty)$.
The regularity properties of the matrix-valued functions $L,L^{-1}$
depend on the functions $A$ and $Q$.
We outline sufficient conditions for the regularity assumptions
of Theorem~\ref{thm:main} in two cases,
where the Cholesky decomposition is considered:
\begin{itemize}
\item[i)]
  Compact domain~$\Pi$ (e.g., uniform distribution, beta distribution, etc): 
  if $A,Q \in \czero^{n \times n}$,
  then $L,L^{-1} \in \czero^{n \times n}$ as shown in the
  proof of Lemma~\ref{lemma:smoothness}.
  Compact domains imply $\czero^{n \times n} \subset \linf^{n \times n}$.
  Hence continuity of $A$ and $Q$ is sufficient.
\item[ii)]
  Unbounded domain~$\Pi$ and exponentially decaying probability 
  density function~$\rho$ for $p \rightarrow \infty$
  (e.g., Gaussian distribution, gamma distribution, etc.): if 
  the components of $A,Q$ are (multivariate) polynomials in~$p$,
  then all entries of both $M$ and $M^{-1}$ are rational functions in~$p$. 
  Consequently, these entries exhibit at most polynomial growth for
  $p \rightarrow \infty$.
  The matrices $L,L^{-1}$ inherit this behavior.
  Since $\rho$ decreases exponentially, we obtain
  $A,L,L^{-1} \in \mathscr{L}^q(\Pi,\rho)$ for any $q \in [1,\infty)$.
\end{itemize}

\subsection{Numerical computation of Galerkin projection}
\label{sec:numerical}
In Theorem~\ref{thm:main}, the transformation matrix $T(p) = L(p)^\top$
depends on the parameters~$p \in \Pi$.
The solution of the Lyapunov equations~(\ref{lyapunov-parameter})
and its symmetric decomposition can be computed analytically
only for simple systems.
We require numerical methods for general systems.

\begin{table} 
  \caption{Probability distributions and
    Gaussian quadrature methods.}
  \label{tab:gauss-quad}
  \begin{center}
    \begin{tabular}{cc}
      probability distribution & quadrature rule \\ \hline
      uniform & Gauss-Legendre \\
      Gaussian & Gauss-Hermite \\
      beta & Gauss-Jacobi \\
      gamma & Gauss-Laguerre
    \end{tabular}
  \end{center} 
\end{table}

We compute the Galerkin projection of the
transformed matrix~(\ref{similarity-trafo}) by a quadrature rule
for the weighted integrals~(\ref{expectedvalue}),
where the probability density is the weight function.
A quadrature scheme is defined by its nodes
$\{ p^{(1)} , \ldots , p^{(k)} \} \subset \Pi$ and weights
$\{ w_1 , \ldots , w_k \} \subset \real$.
For example, Gaussian quadrature can be applied for a single
random variable,
see~\cite[p.~171]{stoerbulirsch}. 
Each traditional probability distribution induces a weighted
integral~(\ref{expectedvalue}) and an associated Gaussian quadrature rule.
Table~\ref{tab:gauss-quad} illustrates the most important cases.
Tensor product rules of Gaussian quadrature can be used for
multiple random variables ($\Pi \subseteq \real^q$),
provided that the number~$q$ is not too large.

Using a general quadrature method,
the approximation $\tilde{B} \in \real^{mn \times mn}$ of
$\hat{B} \in \real^{mn \times mn}$ is given by
\begin{equation} \label{quadrature}
  \hat{B}_{ij} = \mathbb{E} \left[ B \Phi_i \Phi_j \right] \approx
  \tilde{B}_{ij} :=
  \sum_{r=1}^k w_r B(p^{(r)}) \Phi_i(p^{(r)}) \Phi_j(p^{(r)})
\end{equation}
for $i,j=1,\ldots,m$.
The exact matrix $\hat{B} + \hat{B}^\top$ is negative definite due
to Theorem~\ref{thm:main}. 
Given a sufficiently accurate quadrature rule,
the approximation $\tilde{B} + \tilde{B}^\top$ is negative definite as well,
because the eigenvalues of a matrix
depend continuously on its entries.
Hence the linear dynamical system $\dot{\tilde{v}} = \tilde{B} \tilde{v}$
inherits the asymptotic stability.

We show a sufficient condition with respect to the magnitude
of the quadrature error.
\begin{theorem} \label{thm:quadrature}
  Let $\hat{B},\tilde{B},\Delta B \in \real^{mn \times mn}$
  and $\hat{B} = \tilde{B} + \Delta B$.
  If $\hat{B} + \hat{B}^\top$ is negative definite and
  $$ \left\| \Delta B \right\|_2 <
  \textstyle{\frac{1}{2}} | \alpha (\hat{B} + \hat{B}^\top) | $$
  with the spectral (matrix) norm $\| \cdot \|_2$ and
  the spectral abscissa~$\alpha$, then
  $\tilde{B} + \tilde{B}^\top$ is also negative definite.
  Hence the matrix~$\tilde{B}$ is stable.
\end{theorem}

\clearpage

\underline{Proof:} \nopagebreak

The eigenvalues of $\hat{B} + \hat{B}^\top$ are
$\lambda_1 \le \lambda_2 \le \cdots \le \lambda_{mn}<0$.
It holds that $\alpha (\hat{B} + \hat{B}^\top) = \lambda_{mn}$.
The matrix $\hat{B} + \hat{B}^\top$ is symmetric and thus diagonalizable.
An orthonormal basis of eigenvectors exists, which forms a square
matrix of condition number one with respect to the spectral norm.
Let $\mu \in \real$ be an eigenvalue of the symmetric matrix
$\tilde{B} + \tilde{B}^\top$.
The Theorem of Bauer-Fike, see~\cite[p.~357]{golub-loan}, implies
$$ \min_{j=1,\ldots,mn} | \mu - \lambda_i | \le
\left\| \Delta B + \Delta B^\top \right\|_2
\le 2 \| \Delta B \|_2 < | \lambda_{mn} | . $$
The minimum is $| \mu - \lambda_{\ell} |$ with some~$\ell \in \{ 1,\ldots,mn \}$.
It follows that
$ | \mu - \lambda_{\ell} | < - \lambda_{mn} $
and $\mu < - \lambda_{mn} + \lambda_{\ell} \le 0$.
Thus $\tilde{B} + \tilde{B}^\top$ is negative definite.
Lemma~\ref{lemma:stable-matrix} shows that $\tilde{B}$ is a stable matrix.
\hfill $\Box$

\medskip

In Theorem~\ref{thm:quadrature}, the perturbation $\Delta B$ consists of
the quadrature errors concerning~(\ref{quadrature}).
If the quadrature rule is inaccurate, then the matrix
$\tilde{B} + \tilde{B}^\top$ may not be negative definite.
Consequently, the stability can be lost in the approximate
Galerkin projection.

The evaluation of the formula~(\ref{quadrature}) for
all $i,j=1,\ldots,m$ requires mainly to calculate $k$~transformed
matrices~(\ref{similarity-trafo}).
The effort for the evaluation of the basis polynomials is negligible.
We have to solve the
Lyapunov equations~(\ref{lyapunov-parameter}) for $k$~different
realizations of the random parameters.
In addition, we have to compute a symmetric decomposition for each~$M(p^{(r)})$
with $r = 1,\ldots,k$.
Based on the algorithm of Bartels and Stewart~\cite{bartels-stewart},
numerical techniques were derived to compute the Cholesky factor 
without having to compute $M$ first,
see~\cite{hammarling,penzl}.
These direct linear algebra methods exhibit a computational
effort of $O(n^3)$ operations.
Thus our total computational work becomes $O(kn^3)$.
This effort is acceptable in the case of moderate dimensions~$n$.
Typically, we need larger numbers of nodes for polynomials of higher degrees.
For example, with polynomials $\Phi_i$ of degree $i-1$ in
a single random variable ($q=1$),
the approximation~(\ref{quadrature}) requires a Gaussian quadrature
with $k \ge m+1$ nodes.
If the matrix-valued function $B$ is close to a polynomial of degree~$d$,
then $k \approx m + \frac{d}{2}$ nodes are sufficient.

For comparison, an $LU$-decomposition of a dense Galerkin-projected
matrix $\hat{A}$ or $\hat{B}$ costs $O(m^3n^3)$ operations,
which is required in implicit time integrators.
The matrices are typically sparse for high dimensions~$n$,
where iterative methods have to be used for solving Lyapunov equations.
Furthermore, there is some potential to reduce the computational effort by
numerical methods for parameter-dependent Lyapunov equations,
cf.~\cite{son-stykel}.


\section{Nonlinear dynamical systems}
\label{sec:nonlinear}
We derive a stabilization for nonlinear dynamical systems now.

\subsection{Stationary solutions}
Let a nonlinear autonomous dynamical system
\begin{equation} \label{system-nonlinear}
  \dot{x}(t,p) = f(x(t,p),p)
\end{equation}
be given, including a sufficiently smooth function
$f : D \times \Pi \rightarrow \real^n$ ($D \subseteq \real^n$).
We assume that a family of asymptotically stable stationary solutions
$x^* : \Pi \rightarrow \real^n$ exists, i.e.,
\begin{equation} \label{equilibria}
  f ( x^*(p) , p ) = 0 \qquad \mbox{for all} \;\; p \in \Pi .
\end{equation}
The asymptotic stability means that the Jacobian matrix
$\left. \textstyle \frac{\partial f}{\partial x} \right|_{x = x^*(p)}$
is stable for all~$p \in \Pi$, cf.~\cite[p.~22]{seydel}. 

The Galerkin projection of the nonlinear system~(\ref{system-nonlinear})
reads as
\begin{equation} \label{galerkin-nonlinear}
  \dot{\hat{v}}(t) = F(\hat{v}(t))
\end{equation}
with the right-hand side
\begin{equation} \label{galerkin-nl-rhs}
  F = (F_1^\top,\ldots,F_m^\top)^\top , \qquad
  F_i(\hat{v}) := \mathbb{E} \left[ f
  \left( \sum_{j=1}^m \hat{v}_j \Phi_j (\cdot) , \cdot \right)
  \Phi_i(\cdot) \right] 
\end{equation}
for $\hat{v} = (\hat{v}_1^\top , \ldots , \hat{v}_m^\top)^\top$,
where the expected value is applied component-wise again.
However, the existence of an equilibrium of the larger
system~(\ref{galerkin-nonlinear})
with the right-hand side~(\ref{galerkin-nl-rhs})
is not guaranteed in general.
In~\cite{pulch13}, sufficient conditions are specified, under which there is
a stationary solution of~(\ref{galerkin-nonlinear})
for a sufficiently high polynomial degree in gPC expansions.
Moreover, if~(\ref{galerkin-nonlinear}) has a
stationary solution $\hat{v}^* \in \real^{mn}$, then
the function
\begin{equation} \label{equilibria-appr}
  \bar{x}(p) := \sum_{j=1}^m \hat{v}_j^* \Phi_j(p)
\end{equation}
is not an equilibrium of the original system~(\ref{system-nonlinear})
in general.
Thus the stability of the equilibria $x^*$ satisfying~(\ref{equilibria})
does not imply the stability of an equilibrium $\hat{v}^*$
associated to the dynamical system~(\ref{galerkin-nonlinear}).

\subsection{Stabilization of stationary solutions}
\label{sec:stabilization-nonlinear}
Instead of arguing about the stability of an arbitrary equilibrium, we transform 
the system~(\ref{system-nonlinear}) into
\begin{equation} \label{system-nl-shifted}
  \dot{\tilde{x}}(t,p) =
  \tilde{f}(\tilde{x}(t,p),p)
  \qquad \mbox{with} \qquad
  \tilde{f}(x,p) := f(x+x^*(p),p)
\end{equation}
using the family $x^*$ of stationary solutions.
It follows that zero represents an asymptotically stable stationary solution
of the transformed system~(\ref{system-nl-shifted}) for all~$p \in \Pi$.
Let
\begin{equation} \label{galerkin-nl-rhs2}
  \tilde{F} = (\tilde{F}_1^\top,\ldots,\tilde{F}_m^\top)^\top , \qquad
  \tilde{F}_i(\tilde{v}) := \mathbb{E} \left[ \tilde{f}
  \left( \sum_{j=1}^m \tilde{v}_j \Phi_j (\cdot) , \cdot \right)
  \Phi_i(\cdot) \right]
\end{equation}
using the shifted function $\tilde{f}$ from~(\ref{system-nl-shifted})
and $\tilde{v} = (\tilde{v}_1^\top , \ldots , \tilde{v}_m^\top)^\top$.
Now $\tilde{v}^*=0$ is also a stationary solution of the
Galerkin-projected system
\begin{equation} \label{galerkin-nonlinear2}
  \dot{\tilde{v}}(t) = \tilde{F}(\tilde{v}(t))
\end{equation}
with the right-hand side~(\ref{galerkin-nl-rhs2}).
Given a solution of~(\ref{galerkin-nonlinear2}),
an approximation for a solution of the original
system~(\ref{system-nonlinear}) reads as
$$ x(t,p) \approx
\left( \sum_{j=1}^m \tilde{v}_j(t) \Phi_j(p) \right) + x^*(p)
\approx \sum_{j=1}^m \left( v^*_{j} + \tilde{v}_j(t) \right) \Phi_j(p)  $$
assuming a convergent expansion
$$ x^*(p) = \sum_{j=1}^m v^*_{j} \Phi_j(p) $$
of the original stationary solutions.

If the equilibrium $\tilde{v}^*=0$ of the Galerkin-projected
system~(\ref{galerkin-nonlinear2}) is unstable, then a stabilized system
can be constructed as in Section~\ref{sec:galerkin-trafo}.
We apply the transformation to the parameter-dependent matrix
\begin{equation} \label{matrix-A-nonlinear}
  A(p) :=
  \left. \textstyle \frac{\partial f}{\partial x} \right|_{x = x^*(p)} = 
  \left. \textstyle \frac{\partial \tilde{f}}{\partial x} \right|_{x = 0} .
\end{equation}
The stabilized system is given by
\begin{equation} \label{nonlinear-stabilized}
  \dot{y}(t,p) = L(p)^\top \tilde{f} ( L(p)^{-\top} y(t,p) , p ),
\end{equation}
using a symmetric decomposition $M(p) = L(p) L(p)^\top$
of the solution of the Lyapunov equations~(\ref{lyapunov-parameter})
including the parameter-dependent Jacobian matrix~(\ref{matrix-A-nonlinear}).
However, an evaluation of the function~$\tilde{f}$ in the
system~(\ref{nonlinear-stabilized}) requires
the computation of the stationary solution, $x^*$, of~(\ref{system-nonlinear})
for a given parameter value.
Note that the Jacobian matrix of~(\ref{nonlinear-stabilized})
at the stationary solution is
$$ B(p) := L(p)^\top \textstyle
\left. \frac{\partial \tilde{f}}{\partial x} \right|_{x=0}
L(p)^{-\top} . $$
Thus, Section~\ref{sec:galerkin-trafo} yields that the spectrum of the symmetric part of the original Jacobian matrix is changed appropriately to preserve stability.
The Galerkin projection of the dynamical
system~(\ref{nonlinear-stabilized}) results in a larger dynamical
system, whose equilibrium $\tilde{v}^*=0$ is guaranteed to be
asymptotically stable.


\section{Illustrative examples}
\label{sec:example}
We investigate a linear dynamical system as well as a nonlinear dynamical
system with respect to the stability properties
in the Galerkin-projected system. 

\subsection{Linear dynamical system}
We consider the linear dynamical system~(\ref{linear-ode})
including the matrix
\begin{equation} \label{matrix-example}
  A(p) := \textstyle{\frac{1}{100}}
  { \footnotesize
  \begin{pmatrix}
    128 p^2 - 72 p - 32 & 295 p^2 - 199 p + 4 & 165 p^2 - 234 p + 46 \\
    -82 p^2 - 59 p + 270 & -266 p^2 + 144 p - 73 & -147 p^2 - 210 p + 286 \\
    70 p^2 + 296 p - 80 & 43 p^2 + 96 p + 8 & 15 p^2 + 146 p - 251 \\
  \end{pmatrix}
  }
\end{equation}
with a real parameter~$p$.
The eigenvalues of this matrix have a negative real part for
all $p \in [-1,1]$. 
Thus the matrix~(\ref{matrix-example}) is stable in view of
Definition~\ref{def:stable-matrix}.

\begin{figure}
\begin{center}
\includegraphics[width=10cm]{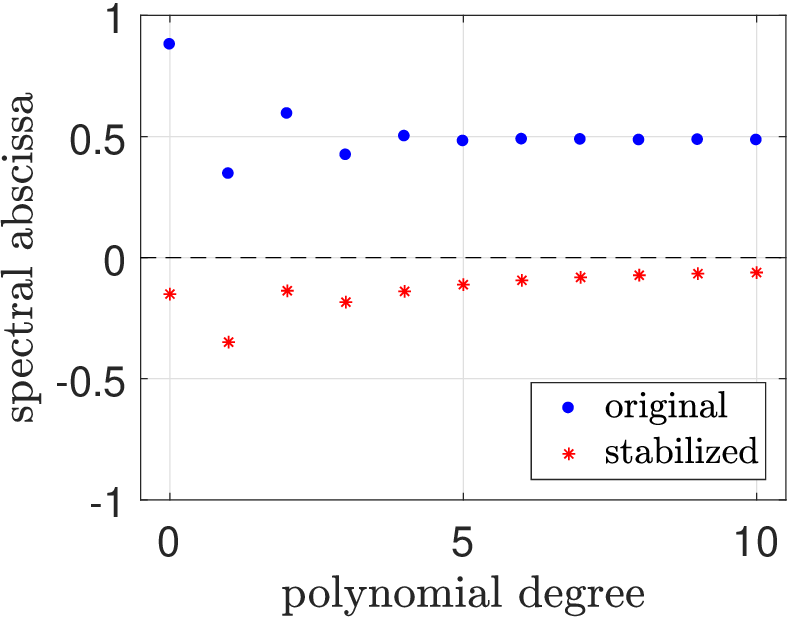} 
\end{center}
\caption{Spectral abscissae in linear dynamical systems,
  including matrix (\ref{matrix-example}), from stochastic Galerkin method
  for different polynomial degrees, using a uniform distribution
  in $[-1,1]$.}
\label{fig:abscissa}
\end{figure}

In the stochastic model, we assume a uniform distribution
for $p \in [-1,1]$.
The expansion~(\ref{expansion}) includes the Legendre polynomials
up to degree~$d = m-1$.
We use Gauss-Legendre quadrature 
to compute the
matrices in the linear dynamical systems~(\ref{system-galerkin})
of the stochastic Galerkin method. This computation is exact, except for round-off errors, because the entries of the matrix~(\ref{matrix-example}) represent
polynomials in~$p$.
However, the Galerkin projection always generates an unstable system.
Figure~\ref{fig:abscissa} illustrates the spectral abscissae
of the matrices~$\hat{A}$ for $d=0,1,\ldots,10$.

Now we use the transformation from Section~\ref{sec:galerkin-trafo}
to obtain a stable system.
In the Lyapunov equation~(\ref{lyapunov-parameter}),
we choose the constant matrix $Q = I$, the identity matrix
in $\real^{3 \times 3}$, which is obviously symmetric and positive definite.
The unique solution $M(p)$ of the Lyapunov equation has entries,
which represent rational functions in the variable~$p$
with numerator/denominator polynomials of degrees up to ten.
The Cholesky algorithm yields the decomposition $M(p) = L(p) L(p)^\top$
with the unique factor.
Hence the transformed matrix~(\ref{B-trafo}) can be computed
point-wise for $p \in [-1,1]$.

The Galerkin projection of the matrix $B(p)$ in the transformed
system~(\ref{system-transformed}) is computed numerically by a
Gauss-Legendre quadrature with 20 nodes.
Thus the matrix~(\ref{B-trafo}) is evaluated at each node of
the quadrature.
Figure~\ref{fig:abscissa} shows the spectral abscissae of the Galerkin-projected matrix for polynomial degrees $d=0,1,\ldots,10$.
We recognize that all spectral abscissae are strictly negative,
which confirms the asymptotic stability.

\begin{figure}
  \begin{center}
  \includegraphics[width=6.5cm]{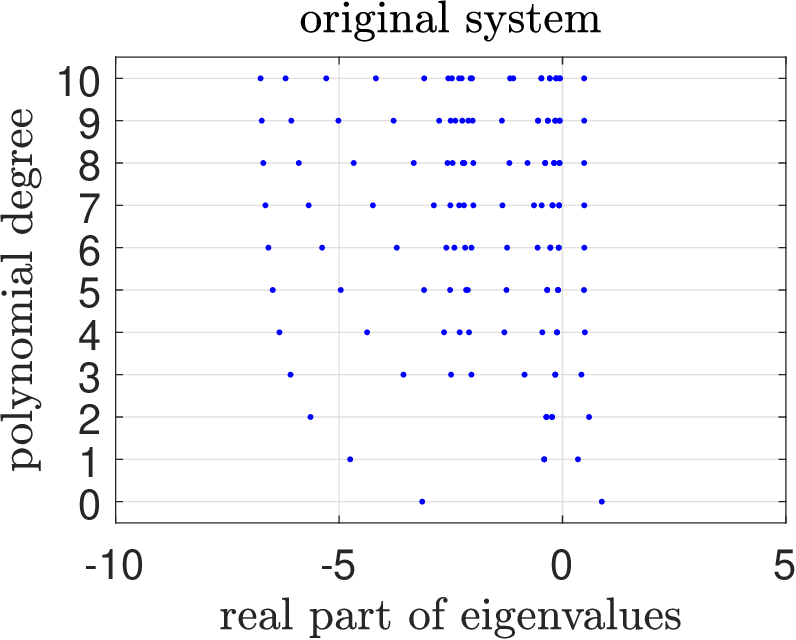}
  \hspace{5mm}
  \includegraphics[width=6.5cm]{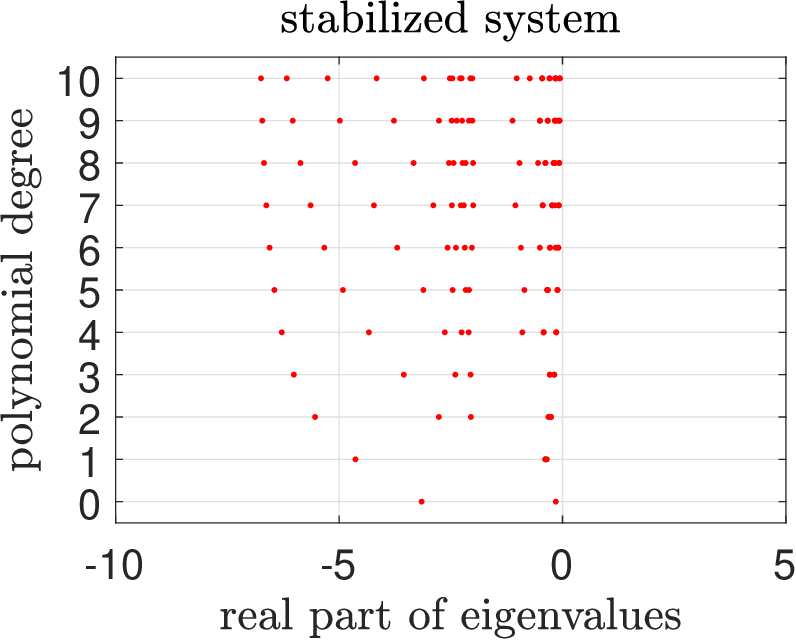}
  \end{center}
  \caption{Eigenvalues of the original system,
    including matrix (\ref{matrix-example}), (left) and the
    stabilized system (right) for different degrees of
    the polynomial expansion with uniform distribution in $[-1,1]$.}
\label{fig:eigenvalues}
\end{figure}

Furthermore, we examine the behavior of all eigenvalues in
the Galerkin-projected systems.
Figure~\ref{fig:eigenvalues} depicts the real part of the eigenvalues
for both the original system and the stabilized system for
different polynomial degrees.
Since complex conjugate eigenvalues arise,
some real parts coincide for each matrix.
We observe that the eigenvalues behave similar and thus
just the stabilization represents the crucial difference. 

\begin{figure} 
\begin{center}
\includegraphics[width=10cm]{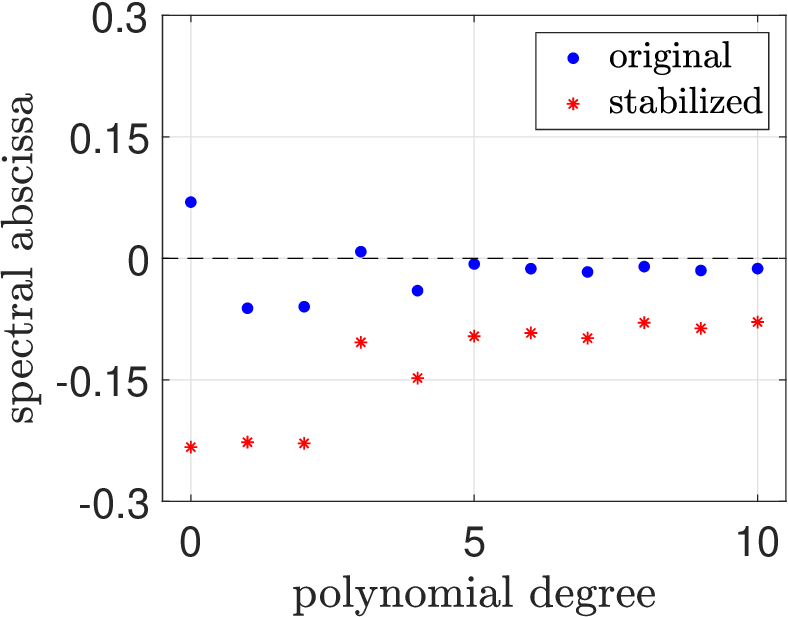} 
\end{center}
\caption{Spectral abscissae in linear dynamical systems,
  including matrix (\ref{matrix-example}), from stochastic Galerkin method
  for different polynomial degrees, using a uniform distribution
  in $[-\frac{2}{5},\frac{2}{5}]$.}
\label{fig:abscissa2}
\end{figure}

We repeat the numerical computations employing a uniform distribution
in the smaller interval $p \in [-\frac{2}{5},\frac{2}{5}]$.
Figure~\ref{fig:abscissa2} illustrates the spectral abscissae of
the Galerkin-projected matrices for different polynomial degrees.
Just two out of ten original systems become unstable now.
Although the other spectral abscissae of the original system are
already negative, the transformed systems exhibit a more negative
spectral abscissa.
In~\cite{pulch19}, a similar stabilization technique was applied in
the context of model order reduction, where this effect becomes
more pronounced in an example.

\begin{figure}
\begin{center}
\includegraphics[width=10cm]{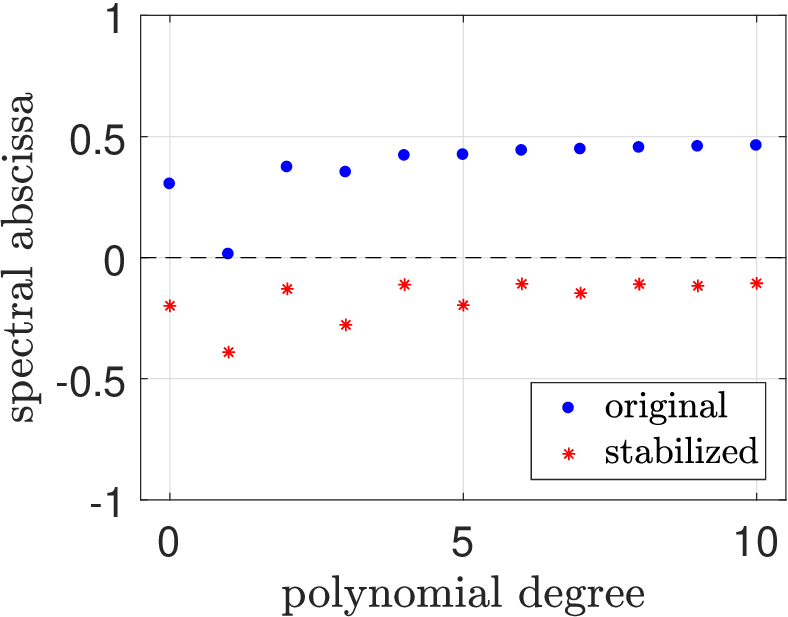} 
\end{center}
\caption{Spectral abscissae in linear dynamical systems, including matrix (\ref{matrix-example}), after the Galerkin projection for different polynomial degrees, using a beta distribution in $[-1,1]$.}
\label{fig:abscissa_beta}
\end{figure}

Alternatively, we choose a beta distribution in the stochastic modeling.
The probability density function reads as
$$ \rho(p) = c \, (1-p)^\alpha (1+p)^\beta \qquad \mbox{for}\;\; p \in [-1,1] $$
with a constant $c>0$ for standardization.
We select $\alpha = 3$, $\beta = 2$.
Jacobi polynomials yield an orthogonal basis.
Again the stochastic Galerkin method results in unstable
systems for all polynomial degrees.
The Galerkin projection of the matrices of the original system and the
stabilized system are computed by Gauss-Jacobi quadrature with 20 nodes. 
The spectral abscissae of the Galerkin-projected matrices 
for the polynomial degrees $d=0,1,\ldots,10$ are depicted in
Figure~\ref{fig:abscissa_beta}.
The transformation yields again a stabilization of the critical systems.

\subsection{Nonlinear dynamical system}
We now consider a two-dimensional dynamical system~(\ref{system-nonlinear}) with
the quadratic right-hand side
\begin{equation} \label{nonlinear-example}
  \mbox{} \hspace{-1.5mm} f(x,p) = \hspace{-1.5mm} \footnotesize
  \begin{pmatrix}
    x_1^2 + (- 35 p - 2 \sin(p) - 13 p^2 - 97) x_1 - 2 x_2^2 +
    (4 \cos(p) - 77 p - 33 p^2 + 23) x_2 \\
    + \sin^2(p) - 2 \cos^2(p) +
    \cos(p) (33 p^2 + 77 p - 23) + \sin(p) (13 p^2 + 35 p + 97) \\[2ex]
    4 x_1^2 + (85 p - 8 \sin(p) + 51 p^2 - 54) x_1 - x_2^2 +
    (2 \cos(p) - \frac{1}{10} p + 67 p^2 - 24) x_2 \\
    + 4 \sin^2(p) - \cos^2(p) +
    \cos(p) (- 67 p^2 + \frac{1}{10} p + 24)
    - \sin(p) (51 p^2 + 85 p - 54) \\
  \end{pmatrix}
\end{equation}
and a real parameter~$p$.
This system is chosen such that
\begin{equation} \label{stationary-example}
  x^*(p) = \begin{pmatrix} \sin (p) \\ \cos (p) \\ \end{pmatrix}
\end{equation}
is a stationary solution for all $p \in \real$.
Numerical computations confirm that these equilibria are
asymptotically stable for all $p \in [-1,1]$.
Since the stationary solution~(\ref{stationary-example}) includes
trigonometric functions, an exact representation in terms of polynomials in the variable~$p$ is not feasible.

Again, we assume the parameter to be a uniformly distributed random variable
with the range $[-1,1]$.
Consequently, the gPC expansion uses the Legendre polynomials.
The stochastic Galerkin method yields the nonlinear dynamical
system~(\ref{galerkin-nonlinear}).
In the right-hand side~(\ref{galerkin-nl-rhs}),
we evaluate the probabilistic integrals (expected values) approximately
by a Gauss-Legendre quadrature using 20 nodes.

\begin{figure}
  \begin{center}  
  \includegraphics[width=6.5cm]{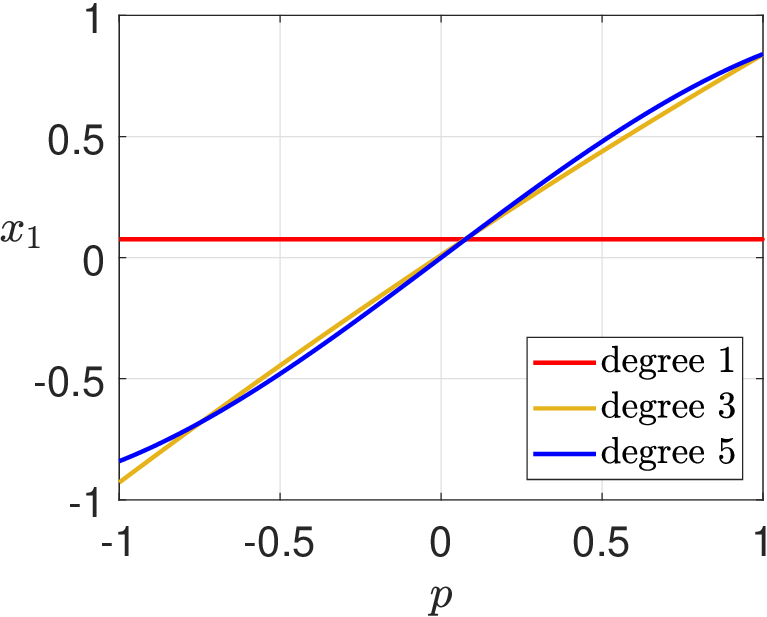}
  \hspace{5mm}
  \includegraphics[width=6.5cm]{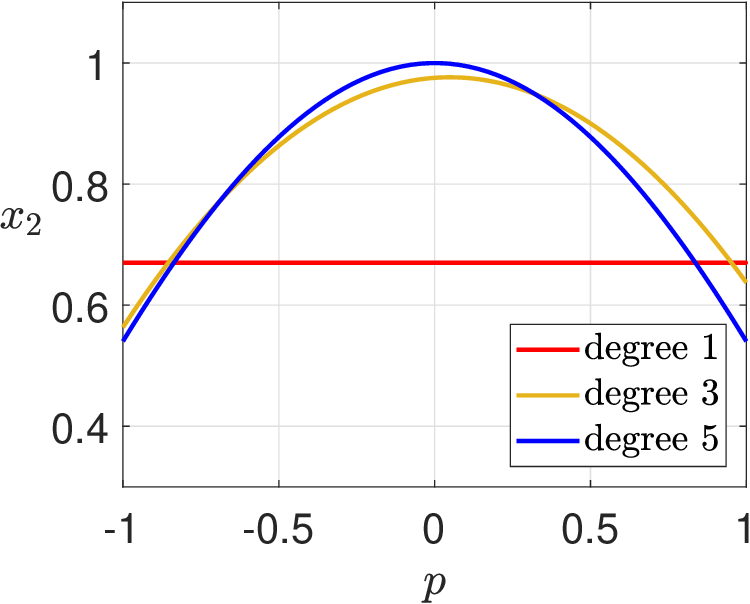}
  \end{center}
  \caption{
    Approximations associated with stationary solutions of the
    Galerkin-projected system, using the right-hand side function (\ref{nonlinear-example}), 
    for different polynomial degrees in quadratic problem.}
\label{fig:nonlinear-equilibria}
\end{figure}

\begin{figure}
  \begin{center}  
  \includegraphics[width=6.5cm]{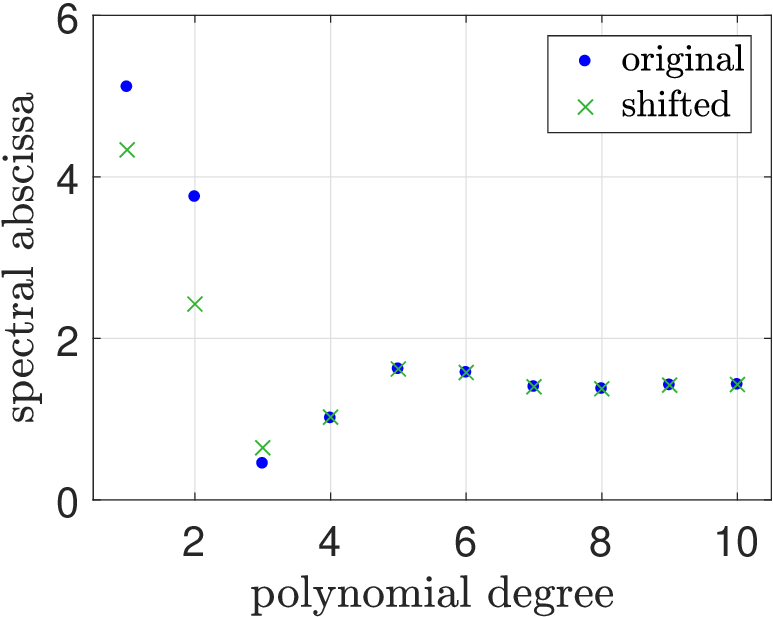}
  \hspace{5mm}
  \includegraphics[width=6.5cm]{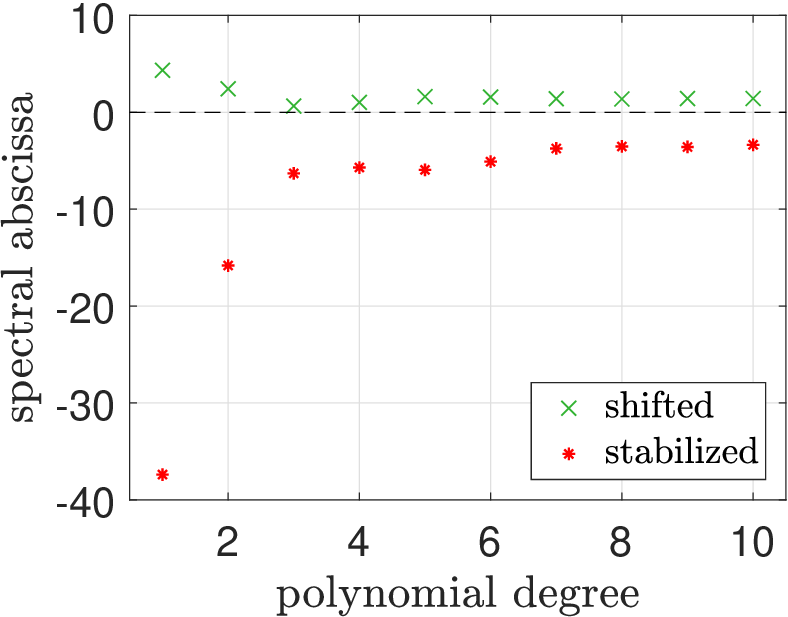}
  \end{center}
  \caption{Spectral abscissae for the Jacobian matrices associated
    with the stationary solutions of the Galerkin projections
    for the original system, the shifted system and the stabilized system 
    with right-hand side function (\ref{nonlinear-example}) in the quadratic problem.}
\label{fig:nonlinear-abscissa}
\end{figure}

Numerical computations show that the Galerkin-projected
system~(\ref{galerkin-nonlinear}) exhibits stationary solutions
for all degrees $d=1,2,\ldots,10$.
Therein, the respective nonlinear systems $F(\hat{v}^*) = 0$ are solved
successfully by Newton iterations.
The corresponding stationary solutions yield the functions~(\ref{equilibria-appr}),
which represent approximations of the equilibrium in the original dynamical system.
Figure~\ref{fig:nonlinear-equilibria} illustrates the
functions~(\ref{equilibria-appr}) for the polynomial degrees $d=1,3,5$.
The approximations converge rapidly to sine and cosine,
respectively, because these trigonometric terms are analytic functions
in the variable~$p$.
However, the stationary solutions of the Galerkin-projected
systems~(\ref{galerkin-nonlinear}) for all polynomial degrees $d$ are 
unstable, which can be seen by the spectral abscissae of the Jacobian matrices
$\left. \frac{\partial F}{\partial \hat{v}} \right|_{\hat{v}=\hat{v}^*}$
in Figure~\ref{fig:nonlinear-abscissa} (left).

To stabilize the computation, we change to the shifted nonlinear dynamical
system~(\ref{system-nl-shifted}) and its Galerkin-projected
system~(\ref{galerkin-nonlinear2}).
The stationary solutions $\tilde{v}^* = 0$ are still unstable
for all polynomial degrees,
which is illustrated by the spectral abscissae of the associated
Jacobian matrices in Figure~\ref{fig:nonlinear-abscissa} (left). 
The values for the original Galerkin system and the novel
Galerkin system become closer for higher polynomial degrees in agreement to
the convergence results in~\cite{pulch13}.

Now we apply the stabilization technique from
Section~\ref{sec:stabilization-nonlinear}.
In the Lyapunov equations~(\ref{lyapunov-parameter}),
the matrix $Q=I$, the identity matrix in $\real^{2 \times 2}$, is selected.
The Cholesky algorithm yields a decomposition of the solutions.
This procedure has to be done for each node of the Gauss-Legendre quadrature.
The stochastic Galerkin method projects the
transformed system~(\ref{nonlinear-stabilized}) to a
larger system~(\ref{galerkin-nonlinear2}).
Figure~\ref{fig:nonlinear-abscissa} (right) shows the spectral abscissae
of the Jacobian matrices for the stationary solution zero
for different polynomial degrees $d=1,2,\ldots,10$.
It follows that the equilibria are asymptotically stable now.

\begin{figure}
  \begin{center}
  \includegraphics[width=8cm]{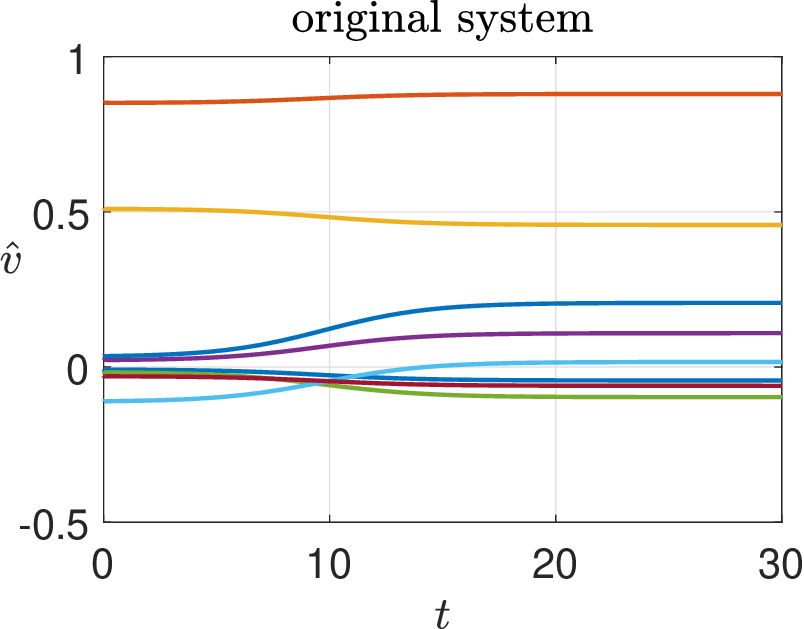}

  \bigskip

  \includegraphics[width=8cm]{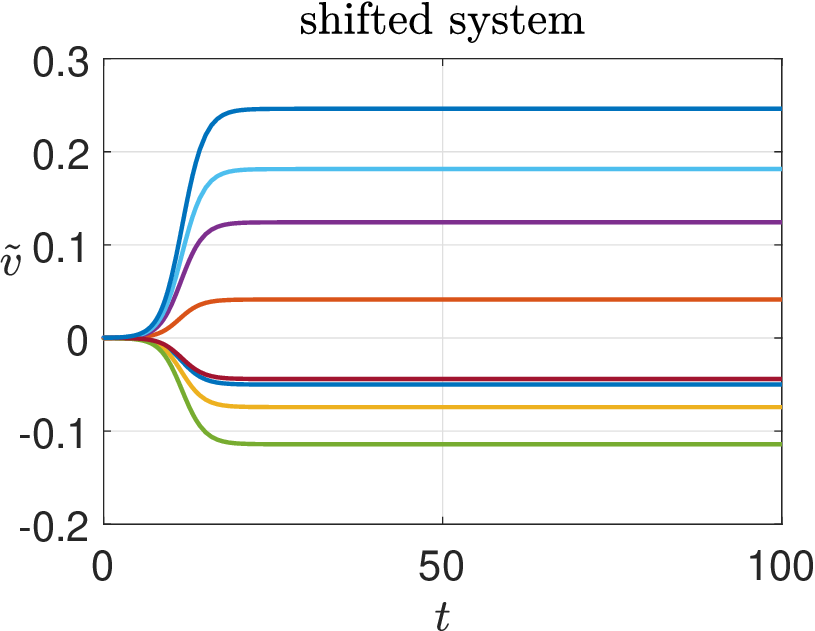}

  \bigskip
  
  \includegraphics[width=8cm]{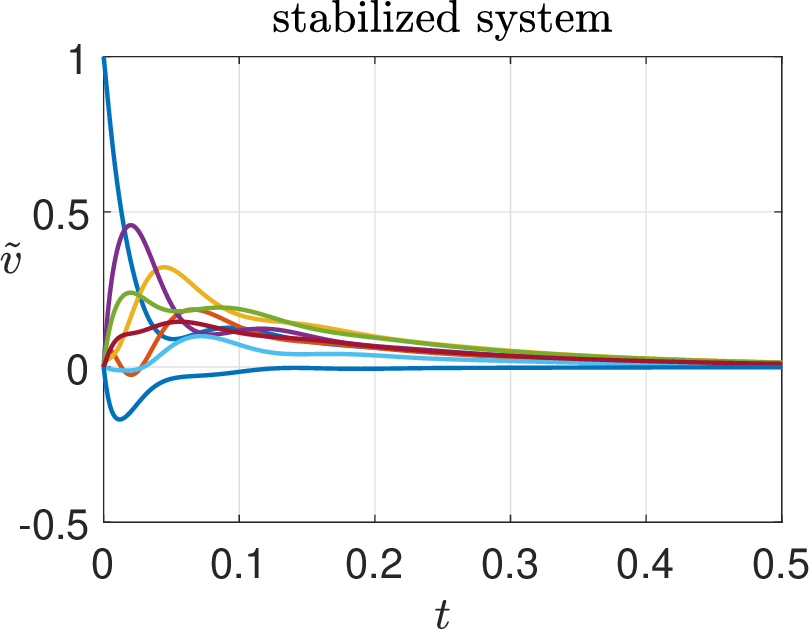}
  \end{center}
  \caption{Solutions of initial value problems for the Galerkin projections
    of original system, shifted system and stabilized system
    in quadratic test example.}
\label{fig:nonlinear-ivp}
\end{figure}

Finally, we illustrate solutions of initial value problems computed using the stochastic Galerkin method for original, shifted and stabilized right-hand side function.
We choose the polynomial degree equal to three,
which results in eight coefficient functions.
The trapezoidal rule yields the numerical solutions.
Firstly, the original Galerkin-projected system~(\ref{galerkin-nonlinear})
is solved, whose initial values are selected close to its stationary solution.
Figure~\ref{fig:nonlinear-ivp} (top) illustrates the numerical solution.
The trajectories are nearly constant at the beginning, say $t \in [0,1]$.
Later the solution changes from the unstable equilibrium
to some stable equilibrium.
Secondly, we solve the Galerkin-projected system~(\ref{galerkin-nonlinear})
with the shifted right-hand side function~(\ref{system-nl-shifted}) using initial values
$\tilde{v}(0) = (10^{-3},0,\ldots,0)^\top$ close to the
equilibrium $\tilde{v}^*=0$.
The same behavior appears like before, as depicted in
Figure~\ref{fig:nonlinear-ivp} (center).
Thirdly, the Galerkin projection of the stabilized
system~(\ref{nonlinear-stabilized}) is solved,
where the initial values are set to $\tilde{v}(0) = (1,0,\ldots,0)^\top$.
Figure~\ref{fig:nonlinear-ivp} (bottom) shows the numerical solution.
Now the trajectories tend to the stable stationary solution
$\tilde{v}^*=0$.

\section{Conclusions}
A basis transformation was constructed for linear dynamical systems
including random variables.
We proved that the stability properties are preserved in
a Galerkin projection of the transformed system.
The transformation matrix follows from a symmetric decomposition of
a solution of a Lyapunov equation.
We showed that the Cholesky factorization retains the smoothness of
involved functions, while the computational effort is low
in comparison to eigenvalue/eigenvector decompositions.
Moreover, the transformation can be applied to guarantee the stability
properties for stationary solutions in nonlinear dynamical systems.
We performed numerical computations for test examples.
The results demonstrate that the combination of the stochastic
Galerkin method and the basis transformation yields an efficient
numerical technique to preserve stability of the original system
under the Galerkin projection.


\end{document}